\documentclass[a4paper]{amsart}
\usepackage{amsfonts,amssymb,amscd,amsmath,amsthm,txfonts}
\usepackage{accents}
\usepackage{txfonts} 
\usepackage{stmaryrd}
\usepackage[all]{xy} 
\usepackage[dvips]{graphicx}
\usepackage{color}

\renewcommand{\restriction}{{\upharpoonright}}
\newcommand{\myeinr}{\phantom{xx}}

\newcommand{\myC}{\mathcal{C}}

\DeclareMathOperator{\mR}{R}

\newtheorem{Thm}{Theorem}[section]
\newtheorem{Lem}[Thm]{Lemma}
\newtheorem{ILem}[Thm]{Induction Lemma}
\newtheorem{Fact}[Thm]{Fact}
\newtheorem{Facts}[Thm]{Facts}
\theoremstyle{definition}
\newtheorem{Def}[Thm]{Definition}
\newtheorem{Asm}[Thm]{Assumption}
\newtheorem{Exm}[Thm]{Example}

\newcommand{\fnto}{\ensuremath{\rightarrow}}  
\DeclareMathOperator{\dom}{dom}               
\newcommand{\card}[1]{{\ensuremath{\rvert#1\lvert}} }     

\newcommand{\om}[1]{\ensuremath{{\omega_{#1}}} }          

\newcommand{\esm}{\ensuremath{\prec}}  

\newcommand{\ho}{\ensuremath{^{\omega}}}                  

\newcommand{\THEN}{\ensuremath{\,\rightarrow\,}}

\newcommand{\DEFEQ}{\ensuremath{\coloneqq}}

\newcommand{\forc}{\ensuremath{\Vdash}}

\newcommand{\n}[1]{\underaccent{\tilde}{#1}}
\newcommand{\an}[1]{\underaccent{\dot}{#1}}

\begin{document}

\date{\today}


\title[New reals]{New reals: Can live with them, can live without them}

\author[Goldstern]{Martin Goldstern}
\address{
  Institut f\"ur Diskrete Mathematik und Geometrie,
  Technische Universit\"at Wien\\
  Wiedner Hauptstra\ss e 8--10/104\\
  1040 Wien, Austria}
\urladdr{http://www.tuwien.ac.at/goldstern}

\author[Kellner]{Jakob Kellner$^\dag$}
\thanks{$^\dag$Partially supported by FWF Austrian Science Fund grant P17627-N12.}
\address{
  Kurt G\"odel Research Center for Mathematical Logic
  at the University of Vienna\\
  W\"ahringer Stra\ss e 25\\
  1090 Wien, Austria}
\urladdr{http://www.logic.univie.ac.at/~kellner}

\begin{abstract}
We give a self-contained proof of the 
preservation theorem for proper countable support iterations 
known as ``tools-preservation'', ``Case~A'' or
``first preservation theorem'' in the literature.
We do not assume that the forcings add reals.\\
\end{abstract}

\maketitle

\section{Introduction}

Judah and Shelah~\cite{JdSh:308} proved that countable support iterations of
proper\footnote{$P$
  is proper if for all countable
  elementary submodels $N\esm H(\chi)$ containing $P$ ($\chi$ a big regular
  cardinal) and all $p\in P\cap N$ there is a $q\leq p$ which forces that $G_P$ is
  $N$-generic (i.e. $G_P\cap D\cap N\neq \emptyset$ for all dense subsets $D\in
  N$). Such a $q$ is called $N$-generic.}
forcings preserve the $\omega\ho$-bounding property (see~\ref{exm:oob} here).
In his book {\em Proper and Improper Forcing\/}~\cite[XVIII \S3]{Sh:f} Shelah
gave several cases of
general preservation theorems for proper countable support iterations (the
proofs tend to be hard to digest, though).  In this paper we deal with
``Case~A''.  

A simplified version of this case appeared in Section 5 of the first 
author's {\em
Tools for your forcing constructions\/}~\cite{tools}. This version uses the
additional requirement that every iterand adds a new real.
Note that this requirement is met in most applications, but the case of
forcings ``not adding reals'' has important applications as well (and note that
not adding
reals is generally not preserved under proper countable support iterations).

A proof of the iteration theorem {\em without\/} this additional requirement
appeared in~\cite{repicky:toolspres} and was copied into {\em Set Theory of the
Reals\/}~\cite{tomekbook} (as ``first preservation theorem'' 6.1.B), but
Schlindwein pointed out a problem in this proof.\footnote{In~\cite{MR2029602},
where Schlindwein gave a proof for the special case of
$\omega^\omega$-bounding, following \cite[VI]{Sh:f}. However he later detected
another problem in his own proof [C. Schlindwein, personal communication, April
2005] and is preparing a new version \cite{CS0505645}.}
In this paper, we
generalize the proof of~\cite{tools}.

We thank Chaz Schlindwein for finding the problems in the existing proofs and
bringing them to our attention.

\section{The Theorem}

Fix a sequence of increasing arithmetical two-place relations 
$(\mR_j)_{j\in\omega}$ on
$\omega\ho$.  Let $\mR$ be the union of the $\mR_j$. Assume
\begin{itemize}
  \item
    $\myC\DEFEQ{}\{f\in \omega\ho:\, f \mR  \eta\text{ for some }
    \eta\in \omega^\omega\}$ is closed,
  \item
    $\{f\in\omega\ho:\, f \,\mR_j\,\eta \}$ is closed for all $j\in\omega$,
    $\eta\in\omega\ho$, and
  \item
    for every countable $N$ there is an $\eta$ such that $ f \mR \eta$ for all
    $f\in N\cap \myC$\\
    (in this case we say ``$\eta$ covers $N$'').
\end{itemize}

\begin{Def}\label{def:martinpreserving}
  Let $P$ be a forcing notion, $p\in P$.
  \begin{itemize}
    \item $\bar{f}^*\DEFEQ{}(f^*_1,\dots,f^*_k)$ is a $P$-interpretation of
      $\bar{\n{f}}\DEFEQ{}(\n{f}_1,\dots,\n{f}_k)$ under $p$, if
      $f^*_i\in\omega\ho$,
      $\n{f}_i$ is a $P$-name for an element of $\myC$, and
      there is an decreasing chain $p\geq p^0\geq p^1\geq \dots$
      of conditions in $P$ such that
      $p^i$ forces $\n{f}_1\restriction i=f^*_1\restriction i\, \&\dots\&\,
      \n{f}_k\restriction i=f^*_k\restriction i$.
    \item A forcing notion $P$ is weakly preserving, if for all\\
      \myeinr $N\esm H(\chi)$ countable, $\eta$ covering $N$, $p\in N$,\\
      there is an $N$-generic $q\leq p$ which forces that $\eta$ covers $N[G_P]$.
    \item A forcing notion $P$ is preserving, if for all\\
      \myeinr $N\esm H(\chi)$ countable, $\eta$ covering $N$, $p\in N$, and\\
      \myeinr $\bar{f}^*,\bar{\n f}\in N$ such that $\bar{f}^*$
        is a $P$-interpretation of $\bar{\n f}$ under $p$,\\ 
      there is an $N$-generic $q\leq p$ which forces that
      $\eta$ covers $N[G_P]$ and moreover that\\
      \myeinr $f^*_i \mR_j \eta$ implies $\n{f}_i \mR_{j} \eta$ for all $i\leq k$, 
              $j\in\omega$.
    \item A forcing notion $P$ is densely preserving if
       there is a dense subforcing $Q\subseteq P$ which is preserving.
  \end{itemize}
\end{Def}

Note that if $\bar{f}^*$ is an interpretation, then $f^*_l\in \myC$
 (since $\myC$ is closed).

The simplest example is that of $\omega\ho$-bounding:
\begin{Exm}\label{exm:oob}
  Set $f \mR_n \eta$ if $f(m)<\eta(m)$ for all $m>n$.
  So $\myC=\omega\ho$, and 
  $f \mR \eta$ if there is an $n$ such that $f(m)<\eta(m)$ for all $m>n$.
  To cover a family of functions means to dominate 
  it. $P$ is weakly preserving iff $P$ is $\omega\ho$-bounding.\footnote{
  $P$ is $\omega\ho$-bounding if
  for all $P$-names $\n f\in\omega\ho$ and $p\in P$ there is a $q\leq p$ 
  and $g\in\omega\ho$ such that $q\forc \n f(m)<g(m)$ for all $m$. 
  So if $P$ is $\omega\ho$-bounding, $\eta$ covers $N$, $\n f\in N$ and
  $G$ is $N$-generic, then $\n f[G]$ is dominated
  by some $g\in N$ and therefore by $\eta$. If on the other hand 
  $P$ is weakly preserving, $\n f$ a $P$-name and  $p\in P$, then
  there is a $N\esm H(\chi)$ containing $p$ and $\n f$. Pick an $\eta\in V$ covering
  $N$. So if $q\leq p$ is as in the definition of weakly preserving,
  then $q$ forces that $\eta$ dominates $\n f$.}
\end{Exm}
This example is typical in the sense that often $\mR$ describes
a covering property of the pair  $(V, V[G])$.

The property ``weakly preserving'' is invariant under equivalent forcings. I.e.
if $P$ forces that there is a $Q$\nobreakdash-generic filter over $V$ and $Q$ forces the
same for $P$, then $Q$ is weakly preserving iff $P$ is weakly
preserving.\footnote{This is analogous (and can be shown analogously)
to the following fact: $P$ is proper (i.e.\ proper for 
all $N\esm H(\chi)$) iff $P$  is proper for all $N\esm H(\chi)$
containing some fixed $x\in H(\chi)$.}
The notion ``preserving'' however does not seem to be invariant.\footnote{
The reason is that the notion of interpretation is not invariant. Given a
forcing $P$ and an interpretation $f^*$ of a function $\n f\notin V$,
we can find a dense subforcing $P'\subset P$ such that for every condition $p'$
of $P'$ there is a $n(p')$ such that $p'$ forces that $f^*(n(p'))\neq \n
f(n(p'))$ (here we identify the $P$-name $\n f$ with the 
equivalent $P'$-name).  So $f^*$ cannot be a $P'$-interpretation of $\n f$.}
It even seems
that ``densely preserving'' does not imply ``preserving''. (Although we do not
have an example. It is not important after all.) 
One direction however is clear:
\begin{Fact}\label{fact:densesubpreservingagain}
  If $P$ is preserving and $Q\subseteq P$ is dense, then $Q$ is preserving.
\end{Fact}

For some instances of $\mR$, weakly preserving is equivalent to preserving.
Most notably this is the case for $\omega\ho$-bounding (see~\cite[6.5]{tools}).

For other instances of  $\mR$ (e.g. Lebesgue positivity, cf.~\cite{KrSh:828})
``$P$ is preserving'' is equivalent to some other property which is invariant under
equivalent forcings.

We will show that densely preserving is preserved under proper
countable support iterations.  This is our version of the theorem known as
``tools preservation''~\cite[Sec. 5]{tools}, ``Case~A''~\cite[XVIII \S3]{Sh:f}
or the
``first preservation theorem''~\cite[6.1.B]{tomekbook}:
\begin{Thm}\label{thm:martin}
  Assume ${(P^0_i,\n{Q}^0_i)}_{i<\epsilon}$ is a countable support iteration
  of proper, densely preserving forcings. Then $P^0_\epsilon$ is
  densely preserving.
\end{Thm}
\newpage
\section{An outline of the proof}

In this section, we describe the ideas used in the proof, without 
being too rigorous.

\subsection*{(A) Use names}

How can we show that the countable support limit of proper 
forcings is proper?

We have a countable support iteration ${(P_\alpha, \n
Q_\alpha)}_{\alpha<\epsilon}$ of proper forcings ($\epsilon$ limit), $N\esm
H(\chi)$ countable, and $p\in P\cap N$. We want to find a $q_\omega\in
P_\epsilon$ which forces that $G$ is $N$-generic, i.e.\ that $G\cap D\cap
N\neq\emptyset$ for all dense subsets $D\in N$ of $P$.

So we fix an
$\omega$-sequence $0=\alpha_0<\alpha_1<\dots$ cofinal in $\epsilon\cap N$,
and enumerate all dense open sets of $P$ that are in $N$ as
${(D_n)}_{n\in \omega}$.

\begin{figure}[tb]
  \begin{minipage}{0.49\textwidth}
    \begin{center}
      \scalebox{0.4}{\input{tries_first.pstex_t}}
      \caption{\label{fig:try1}}
    \end{center}
  \end{minipage}
  \hfill 
  \begin{minipage}{0.49\textwidth}
    \begin{center}
      \scalebox{0.4}{\input{tries_point.pstex_t}}
      \caption{\label{fig:try2}}
    \end{center}
  \end{minipage}
\end{figure}
One unsuccessful attempt to construct $q_\omega$ could be the one
illustrated in Figure~\ref{fig:try1}:
Set $p_{-1}\DEFEQ p$ and $q_{-1}\DEFEQ \emptyset$.
Given $p_{n-1}\in N$ and $q_{n-1}$, 
choose (in $N$) a $p_n\leq p_{n-1}$ in $D_n\cap N$ and 
 (in $V$) a
$q_n\leq p_n\restriction \alpha_{n+1}$ which extends $q_{n-1}$.
Set $q_\omega\DEFEQ \bigcup q_n$. Then
$q_\omega$ is $N$-generic, since $q_\omega\leq p_n\in D_n\cap N$.
Of course this doesn't work, since we generally cannot
find a $p_n\leq p_{n-1}$ in $D_n$ 
such that $q_{n-1}\leq p_n\restriction \alpha_n$.

What we actually do instead is the following (see Figure~\ref{fig:try2}):
The $p_n$ will be $P_{\alpha_n}$-names, and the $q_n$ are
$P_{\alpha_{n+1}}$-generic over $N$. So
instead of choosing \mbox{$p_n\in P_{\epsilon}$}, we 
choose (in $N$) a $P_{\alpha_n}$-name $\n p_n$ for an element of $P_{\epsilon}$
such that the following is forced by $P_{\alpha_n}$:
\begin{itemize} 
  \item $\n p_n\in D_n$,
  \item $\n p_n\restriction \alpha_n\in G_{\alpha_n}$, and
  \item if $\n p_{n-1}\restriction \alpha_n\in G_{\alpha_n}$,
        then $\n p_n\leq \n p_{n-1}$.
\end{itemize} 
It is clear that we can find such a name. So we first construct all the
$\n p_n$ (each $\n p_n$ is in $N$, but the sequence is not).
Then we construct $q_n\in P_{\alpha_{n+1}}$ satisfying the following:
\begin{itemize} 
  \item $q_{n}$ extends $q_{n-1}$,
  \item $q_{n}$ is $P_{\alpha_{n+1}}$-generic over $N$, and
  \item $q_{n}$ is stronger than
        $\n p_{n}$ on the interval $[\alpha_{n},\alpha_{n+1})$.\footnote{More
     formally (since $\n p_{n}$ is a name): 
     For all $\alpha_n\leq \beta < \alpha_{n+1}$,
      $q_{n}\restriction \beta\forc_\beta p_n\restriction \beta\in G_\beta\ \& \
      q_n(\beta)\leq \n p_n(\beta)$.
  }
\end{itemize} 
So (by induction) $q_n$ forces that $\n p_n\restriction \alpha_{n+1}\in
G_{\alpha_{n+1}}$ and that therefore $\n p_{n+1}\leq \n p_n$. 
So $q_\omega=\bigcup q_n$ 
forces that $\n p_n\restriction \alpha_n\in G_\alpha$ (by
definition of $\n p_n$), that $\n p_n\restriction \alpha_{n+1}\geq 
\n p_{n+1}\restriction \alpha_{n+1}\in G_{\alpha_{n+1}}$ and generally
that $\n p_n\restriction \alpha_{m}\in G_{\alpha_{m}}$ for all $m>n$.
Therefore $q_\omega$ forces that $\n p_n\in G_\epsilon$.
Also, $q_{n-1}$ is $P_{\alpha_n}$-generic over $N$, and the
$P_{\alpha_n}$-name $\n p_n$ is in $N$,
so $q_\omega$ forces that $\n p_n\in N\cap P_{\epsilon}$ and therefore
in $N\cap D_n\cap G_\epsilon$,  i.e.\ that $G_\epsilon$ is $N$-generic.

\subsection*{(B) Interpolate approximations}
First note that for every $P_\epsilon$-name $\n f\in\myC$ and for 
every $p\in P_\epsilon$ we can find an approximation $f^*$ of
$\n f$ under $p$. If additionally $0<\alpha<\epsilon$ and $P_\alpha$
adds a new real $\n r$, then we can choose the witnesses of the approximation
such that $\{p^m\restriction \alpha:\, m\in\omega\}\subseteq P_\alpha$ 
is inconsistent.\footnote{We call a set $A\subseteq P$ inconsistent,
if $P$ forces that not every condition of $A$ is in $G$.} 
(Just let $p^m\restriction \alpha$ decide $\n r(m)$.)

Now
assume that $f^*$ is a $P_\epsilon$-approximation of $\n f$ witnessed by 
$(p^m_0)_{m\in\omega}$ and that 
$\{p^m_0\restriction \alpha:\, m\in\omega\}\subseteq P_\alpha$ is inconsistent.
Then we can define $P_\alpha$-names $(\n p^m_\alpha)_{m\in\omega}$ and $\n
f^{**}$ such that the following is forced by $P_\alpha$ (see Figure~\ref{fig:interpolation}):
\begin{figure}[tb]
  \begin{minipage}{0.35\textwidth}
    \begin{center}
      \scalebox{0.4}{\input{interpolation.pstex_t}}
      \caption{\label{fig:interpolation}}
    \end{center}
  \end{minipage}
  \hfill 
  \begin{minipage}{0.62\textwidth}
    \begin{center}
      \scalebox{0.4}{\input{total.pstex_t}}
      \caption{\label{fig:total}}
    \end{center}
  \end{minipage}
\end{figure}

\begin{itemize}
\item $p^m_0\restriction \alpha\in G_\alpha$ implies $\n p^0_\alpha\leq p^m_0$
  (i.e. $\n p^0_\alpha$ is stronger than the strongest $p^m_0$ whose restriction
  is in $G_\alpha$),
\item $\n f^{**}$ is an approximation of $\n f$ witnessed by 
  $(\n p^m_\alpha)_{m\in\omega}$.
\end{itemize}
Then $(p^m_0\restriction \alpha)_{m\in\omega}$ 
witnesses that $f^*$ approximates $\n f^{**}$:

$p^m_0\restriction \alpha$ forces that
\begin{itemize}
\item
 $\n p^m_\alpha$ forces that $\n f^{**}\restriction m=\n f\restriction m$
  and
\item
  $\n p^m_\alpha\leq p^m_0$ and therefore that
\item $\n p^m_\alpha$ 
  also forces $f^*\restriction m=\n f\restriction m$.
\end{itemize}
So $p^m_0\restriction \alpha \wedge \n p^m_\alpha$ forces 
$\n f^{**}\restriction m=\n f\restriction m=f^*\restriction m$, and
since $\n f^{**}\restriction m$, $f^*\restriction m$ already live in $V[G_\alpha]$,
$\n f^{**}\restriction m=f^*\restriction m$ is already forced by
$p^m_0\restriction \alpha$. 

So we can interpolate (or ``factorize'') the interpretation 
$(f^*,\n f)$ by the ``composition'' of the interpretations
$(f^*,\n f^{**})$ and $(\n f^{**},\n f)$.

\subsection*{(C) Approximate more and more functions better and better}
In addition to all the dense sets $D_n$ of $N$ --- as in (A) --- we also list
all the $P_\epsilon$-names $\n f_n$ in $N$ for elements of $\myC$. We have
to make sure that $q_\omega$ forces that $\n f \mR \eta$.  We assume that
every element of 
$D_n$ decides $\n f_m\restriction n$ for $m\leq n$.

We start with an approximation $f^{*0}_0$ for $\n f_0$ witnessed
by $(p^m_0)_{m\in\omega}$. We assume that 
\mbox{$\{p^m_0\restriction \alpha_1:\, m\in\omega\}$} is inconsistent.
We can 
find (in $N$) $P_{\alpha_1}$ names $(\n p^m_1)_{m\in\omega}$ and 
$\n f^{*1}_0, \n f^{*1}_1$
(see Figure~\ref{fig:total}) 
such that the following is forced:
\begin{itemize}
\item $\n f^{*1}_0, \n f^{*1}_1$ are interpretations of 
  $\n f_0,\n f_1$ witnessed by $(\n p^m_1)_{m\in\omega}$,
\item $p^m_0\in G_{\alpha_1}$ implies $\n p^0_1\leq p^m_0$ 
   (i.e. $\n f^{*1}_0$ interpolates $(f^{*0}_0,\n f_0)$ as in (B)),
\item $\n p^0_1\in D_1$ (in particular, $\n p^0_1$ decides 
       $\n f_0\restriction 1,
       \n f_1\restriction 1$), and
\item we again assume that 
  $\{\n p^m_1\restriction \alpha_2:\, m\in\omega\}$ is inconsistent.
\end{itemize}
Because of the last item, we can iterate this construction.

Now we choose (in $V$) 
a $q_0\in P_{\alpha_1}$ such that $q_0\leq p^0_0\restriction\alpha_1$
and $q_0$ is $P_{\alpha_1}$-generic over $N$ and forces that
$\eta$ covers $N[G_{\alpha_1}]$ and that $f^{*0}_0 \mR_j \eta$
implies $\n f^{*1} \mR_j \eta$ for all $m$. Inductively, we get
a sequence $(q_n)_{n\in\omega}$ such that $q_n\in P_{\alpha_{n+1}}$
extends $q_{n-1}$ and
forces
\begin{itemize}
\item $G_{\alpha_{n+1}}$ is $N$-generic and $\eta$ covers $N[G_{\alpha_{n+1}}]$,
\item $\n f^{*n}_m \mR_j \eta$ implies $\n f^{*n+1}_m \mR_j \eta$ for $m\leq n$ and
	all $j$.
\end {itemize}

Let $q_\omega$ be the union of all $q_n$. Then $q_\omega$ forces the following:
For $m\geq n$, $\n f_n\restriction m=\n f^{*m}_n\restriction m$ (since $\n
p^0_m\in D_m$ decides $\n f_n\restriction m$).  Also, $\n f^{*n}_n \mR_j \eta$ for
some $j\in\omega$ (since $\n f^{*n}_n\in N[G_{\alpha_n}]$ and $\eta$ covers
$N[G_{\alpha_n}]$).  $\n f_n$ is the limit of functions $\n f^{*m}_n$ which all
satisfy $\n f^{*m}_n \mR_j \eta$. Since $\{f\in\omega\ho:\, f \,\mR_j\,\eta \}$ is
closed, $\n f_n \mR_j \eta$. Also, $q_\omega$ is $N$-generic just as in (A).

\subsection*{(D) Decide when we are $\sigma$-complete}

The proof so far relies on the fact that we can always find approximations
whose witnesses are inconsistent.

We already know that this is the case if the iteration between $\alpha_n$ and
$\alpha_{n+1}$ adds a new real.
Actually we just need that the iterands are ``nowhere
$\sigma$-complete'', i.e. that below every $p$ we can find an inconsistent
decreasing sequence.

If no reals are added, it might seem as we do not have anything to do
(since Case~A preservation is vacuous without new reals). 
The problem is that the
countable support iteration of proper forcings which do not add
reals can add a real in the limit.
So it might be that we are unable to use new reals in the intermediate steps (which we want to
construct inconsistent witnesses for approximations), but get
new reals in the limit (which could be a problem for preservation).

On the other extreme, if all iterands are $\sigma$-complete, then the
limit is $\sigma$-complete as well, and therefore adds no reals,
so there is nothing to do.

So what to do?

First note that we can split every forcing in a $\sigma$-complete and a nowhere
$\sigma$-complete part. However, that does not solve our problem,
since we can not split the index set $\epsilon$ of the iteration into
$\epsilon_1,\epsilon_2$ such that $P_\alpha$ forces that $\n Q_\alpha$ is
$\sigma$-complete if $\alpha\in \epsilon_1$ and nowhere $\sigma$-complete
otherwise.

For example, $Q_0$ could add a Cohen real $\n c$, and $\n Q_n$ could be defined
to be $\sigma$-complete iff $\n c(n)=0$.

So we will do the following: Given a condition $p\in P_\epsilon$,
there is a maximal $\gamma\leq \epsilon$ such that 
$P_\alpha$ forces that $\n Q_\alpha$ is $\sigma$-complete (below $p(\alpha)$)
for all $\alpha<\gamma$. So if $\gamma= \epsilon$, then the
rest of the iteration is $\sigma$-complete. If $\gamma< \epsilon$, 
then we strengthen $p$ such that 
$P_\gamma$ forces that $\n Q_\gamma$  is nowhere $\sigma$-complete
(below $p(\gamma)$).

We will only be interested in honest approximations, that is an approximation
witnessed by $(p^m)_{m\in\omega}$ where $p^0$ (and therefore all $p^m$) 
will know
the $\gamma$ where $\n Q_\alpha$ stops to be $\sigma$-complete (in the
way just described).

Since in (C) the conditions $\n p^m_n$ are $P_{\alpha_n}$-names, the 
corresponding $\gamma$ will be a $P_{\alpha_n}$-name as well.
In the iteration at stage $n$, we will have to distinguish three cases:
\begin{itemize}
\item $\{\n p^m_{n-1}\restriction \alpha_n\}$ is inconsistent. Then 
  continue as in (B).
\item The $\gamma$ corresponding to $\n p^0_{n-1}$ is bigger than
  $\alpha_n$ but less than $\epsilon$.
  Then just ``do nothing'', i.e. wait in the iteration
  until $\alpha_m$ is above $\gamma$ and therefore the 
  witnesses are inconsistent.
\item Otherwise, we know that the rest of the iteration is $\sigma$-complete.
\end{itemize}

Again, we do not know from the beginning which case we will use at a given
stage.  In the example above, we will do nothing at stage $n$ iff $\n c(n)=0$
(so it will never happen that the rest of the iteration is $\sigma$-complete).

Also, when we ``do nothing'', we cannot increase the number of functions we
approximate. In (C), the number $k(n)$ of functions which we approximate in step
$n$ was $n+1$ ($\n f^{*n}_0,\dots, f^{*n}_n$ approximates $\n f_0,\dots, f_n$).
So in the proof this number $\n k_n$ will be a $P_{\alpha_n}$-name which is
$\n k_{n-1}$ in case ``do nothing'' and $n+1$ otherwise.

\section{The proof}

\begin{Def}
  Let $Q$ be a forcing, $q\in Q$.
  \begin{itemize}
  \item $q$ is $\sigma$-complete in $Q$, if $Q_q\DEFEQ\{r\in Q:\, r\leq q\}$
    is $\sigma$-complete. In this case we write $q\in Q^\sigma$.
  \item $q$ is nowhere $\sigma$-complete in $Q$
    if there is no $q'\leq_Q q$ such that $q'\in Q^\sigma$.
    In this case we write $q\in Q^{\lnot\sigma}$.
  \item
    $Q$ is decisive if every $q\in Q$ 
    is either $1_Q$ (the weakest element of $Q$)
    or $\sigma$-complete or nowhere $\sigma$-complete.\footnote{Of
    course it is possible to have 
    $1_Q\in Q^\sigma$ or $1_Q\in Q^{\lnot\sigma}$.}
  \end{itemize}
\end{Def}

\begin{Fact}\label{fact:decisiveisdense}
  For every $P$ the set of conditions that are either $\sigma$-complete or
  nowhere $\sigma$-complete is open dense. I.e. for every $P$ there is a dense
  subforcing $Q\subseteq P$ which is decisive.
\end{Fact}

\begin{Fact}\label{denseiswurscht}
  If $(P_\alpha,\n Q_\alpha)_{\alpha<\epsilon}$ is an iteration and 
  $P_\alpha$ forces that $\n Q''_\alpha\subseteq \n Q_\alpha$ is dense
  (for every $\alpha\in\epsilon$),
  then there are an iteration $(P'_\alpha,\n Q'_\alpha)_{\alpha<\epsilon}$ 
  and dense embeddings $\varphi_\alpha:P'_\alpha\fnto P_\alpha$ 
  ($\alpha\leq \epsilon$) such that for $\alpha\leq \beta\leq \epsilon$
  the following holds:
  \begin{itemize}
    \item If $p\in P'_\beta$ then 
      $\varphi_\alpha(p\restriction \alpha)=\varphi_\beta(p)\restriction\alpha$.
    \item In particular $\varphi_\beta$ is an extension of $\varphi_\alpha$.
    \item $P'_\alpha$ forces that $\n Q'_\alpha=\n Q''_\alpha[G_{P_\alpha}]$.\footnote{Where $G_{P_\alpha}\DEFEQ  \{p\in P_\alpha:\, 
         (\exists p'\in  G_{P'_\alpha})\, \varphi_\alpha(p')\leq p\}$ is
the canonic $P_\alpha$-generic filter over $V$.} 
  \end{itemize}
\end{Fact}

Because of \ref{fact:densesubpreservingagain}, \ref{fact:decisiveisdense} and
\ref{denseiswurscht} we can modify the original iteration $(P^0_\alpha,\n
Q^0_\alpha)_{\alpha<\epsilon}$ of Theorem \ref{thm:martin}
to get an iteration $(P_\alpha,\n
Q_\alpha)_{\alpha<\epsilon}$ satisfying $P_\epsilon$ is a dense subforcing of $P^0_\epsilon$ and:

\begin{Asm}\label{Qisstraight}
  $P_\alpha$ forces that $\n Q_\alpha$ is proper,
  decisive and preserving.  
\end{Asm} 

We will show that in this case $P_\epsilon$ is densely
preserving,\footnote{Note that we do not claim that $P_\epsilon$ is
preserving.} so $P^0_\epsilon$ is densely
preserving as well, proving Theorem \ref{thm:martin}.

{\bf From now on we fix} the iteration $(P_\alpha,\n
Q_\alpha)_{\alpha<\epsilon}$ satisfying \ref{Qisstraight}.  We also fix a
regular \mbox{$\chi\gg 2^\card{P_\epsilon}$}, a countable $N\esm H(\chi)$
containing $(P_\alpha,\n Q_\alpha)_{\alpha<\epsilon}$, and an $\eta$
covering $N$. 

\begin{Def}
  We will use the following notation 
  ($\alpha\leq \beta$):
  \begin{itemize}
  \item For $p\in G_\alpha$, $p \forc_\alpha \varphi$ means
    $p \forc_{P_\alpha} \varphi$.
  \item If  $p\in G_\beta$, $r\in P_\alpha$ and $r\leq p\restriction \alpha$, then we can
    define $r\wedge p\in G_\beta$, the weakest condition stronger than $r$ and $p$.
  \item $G_\alpha$ is the $P_\alpha$-generic filter over $V$
    (or its canonical name). So $\forc_\beta G_\alpha=G_\beta\cap P_\alpha$.
    We set $V_\alpha\DEFEQ V[G_\alpha]$.
  \item
    $P_\beta/G_\alpha$ is the $P_\alpha$-name for the forcing consisting of
    those $P_\beta$-conditions $p$ such that $p\restriction\alpha\in G_\alpha$
    (with the same order as $P_\beta$).
  \item In $V_\alpha$: If $p\in P_\beta/G_\alpha$, then
    $p\forc_{(\alpha,\beta)}\varphi$ means
    $p\forc_{P_\beta/G_\alpha}\varphi$. We also say 
    ``$p$ $(\alpha,\beta)$-forces $\varphi$''.
  \end{itemize}
\end{Def}
\begin{Facts}
  Let $0\leq \alpha\leq \beta\leq \epsilon$.
  \begin{itemize}
  \item
    The function $P_\beta\fnto P_\alpha\ast P_\beta/G_\alpha$ 
    defined by $p\mapsto (p\restriction \alpha,p)$ is a dense embedding.
  \item
    If $p_1\in P_\alpha$ and $\n p_2$ is a $P_\alpha$-name for an element 
    of $P_\beta/G_\alpha$, then 
    $p_1\forc_\alpha \n p_2\forc_{(\alpha,\beta)}\varphi$ is equivalent to
    \mbox{$\forc_{\beta} (p_1\in G_\beta\,\&\,\n p_2\in G_\beta)\ \THEN\ \varphi$}.
  \item\label{item:dense}
   If $D$ is an (open) dense subset of $P_\beta$, then
      $D\cap P_\beta/G_\alpha$  is a $P_\alpha$-name for an
      (open) dense subset of $P_\beta/G_\alpha$.
  \end{itemize}
\end{Facts}

If $\n p$ is a $P_\alpha$-name for an element of 
$P_\beta/G_\alpha$, then $\forc_\alpha \n p\forc_{(\alpha,\beta)} \varphi$
does not imply that $\n p[G_\alpha]$ (which is an element of 
$P_\beta$ and therefore of $V$) forces $\varphi$ in $V$ (as element 
of $P_\alpha$). I.e. $V\vDash( \forc_\alpha \n p \forc_{(\alpha,\beta)}\varphi)$
does not imply $\forc_\alpha(V\vDash \n p \forc_{\beta}\varphi)$.

We will use the following straightforward technical facts:
\begin{Lem}\label{lem:quotient}
  Let $0\leq \alpha\leq \gamma\leq\beta\leq \epsilon$. $P_\alpha$ forces:
  \begin{enumerate} 
  \item\label{item:modify1}
     If $p\in P_\beta/G_\alpha$, $q\in P_\gamma/G_\alpha$,
     and $q\forc_{(\alpha,\gamma)} p\restriction \gamma\in G_\gamma$,
     then we can define $p'=q\wedge (p\restriction \beta\setminus\gamma)$ 
     in $P_\beta/G_\alpha$ such that
     $p'\restriction \gamma=q$ and 
     $p'\restriction \xi\forc_{(\alpha,\xi)}p'(\xi)=p(\xi)$
     for $\gamma\leq\xi<\beta$. If 
     $q\leq p\restriction\gamma$ then $q\wedge (p\restriction
     \beta\setminus\gamma)\leq p$,
     and if $p_2\leq p_1$ then $q\wedge (p_2\restriction
     \beta\setminus\gamma)\leq 
     q\wedge (p_1\restriction \beta\setminus\gamma)$.
  \item\label{item:triv3252}
    If $p^0\geq p^1\geq \dots$ is a decreasing sequence in
    $P_\gamma/G_\alpha$, and
    for every $\alpha\leq\zeta< \gamma$ we have
    $p^0\restriction \zeta \forc_{(\alpha,\zeta)}  p^0(\zeta)\in \n
    Q^\sigma_\zeta$,
    then there is a $p^\omega\leq p^0\in  P_\gamma/G_\alpha$
    such that $p^\omega \forc_{(\alpha,\gamma)} p^m\in G_\gamma$ for all $m\in
    \omega$. (Here we actually use that $P_\alpha$ is proper.)
  \end{enumerate}
\end{Lem}

\begin{proof}
  To show (\ref{item:modify1}), set $A\DEFEQ
  \dom(q)\cup(\dom(p)\setminus\alpha)$.
  Note that $A\in V$. Fix a $P_\alpha$-name for $p$.
  Define for $\xi\in A$ (in $V$) $p'(\xi)=q(\xi)$ if $\xi<\gamma$,
  and for $\xi\geq \gamma$ let $p'(\xi)$ be $p(\xi)$ provided that
  $p\restriction\xi\in G_\xi$ ($1_{\n Q_\xi}$ otherwise).

  (\ref{item:triv3252}) is similar:
  There is a $A\in V$ countable in $V$ such that $A\supseteq
  \bigcup_{m\in\omega}
  \dom(p^m)$ (since $P_\alpha$ is proper). Fix a $P_\alpha$-name (in $V$) 
  for the sequence $(p^m)_{m\in\omega}$.
  
  Now define $p^\omega$ in $V$:
  Set $p^\omega \restriction \alpha \DEFEQ p^0\restriction \alpha$.
  For $\alpha\leq \zeta<\gamma$, $\zeta\in A$ define
  $p^\omega(\zeta)\in \an Q_\zeta$ to be a lower
  bound of $\{p^m(\zeta):\, m\in\omega\}$ if such a lower bound exists,
  and $p^0(\zeta)$ otherwise.
\end{proof}

 From now on, to distinguish between $P_\beta$-names and $P_\alpha$-names for
some $\alpha<\beta$, we denote $P_\beta$-names (in $V$ as well as
$P_\alpha$-names
for such names) with a tilde under the symbol (e.g. $\n\tau$) and we denote
$P_\alpha$-names for $V_\alpha$ objects that are not $P_\beta$-names (but could
be $P_\beta$ conditions) with a dot under the symbol (e.g. $\an\tau$). 
In particular we write $(P_\alpha,\an Q_\alpha)_{\alpha<\epsilon}$.

\begin{Def}
  Let $\alpha\leq \beta\leq \epsilon$. Work in $V_\alpha$.
  \begin{itemize}
  \item
    $(p^m)_{m\in\omega}$ is an honest $(\alpha,\gamma,\beta)$-sequence, if
    \begin{itemize}
    \item $p^m\in P_\beta/G_\alpha$, 
    \item $p^{m+1}\leq p^m$,
    \item $\alpha\leq \gamma\leq\beta$,
    \item for all $\alpha \leq \zeta <\gamma$,
      $p^0\restriction\zeta\forc_{(\alpha,\zeta)}p^0(\zeta)\in 
      \an Q_\zeta^{\sigma}$,\footnote{if $\zeta\notin\dom(p)$,
        then $p(\zeta)$ is defined to be $1_{Q_\zeta}$. In this case
        $p(\zeta)\in Q_\zeta^{\sigma}$ means that $Q_\zeta$ 
        is $\sigma$-complete.
        So it is possible that $\gamma\geq \alpha+\om1$,
        this is no contradiction to countable support.} 
    \item 
      $p^m\restriction\gamma=
      p^0\restriction\gamma$ for all $m$.
    \item if $\gamma<\beta$, then $p^0\restriction\gamma$ 
      $(\alpha,\gamma)$-forces that\\
      \myeinr $p^0(\gamma)\in \an Q_{\gamma}^{\lnot\sigma}$, and\\
      \myeinr 
        $\{p^m(\gamma):\,m\in\omega\}\subseteq \an Q_{\gamma}$ is inconsistent.
    \end{itemize}
  \item   
    Let $k$ be a natural number, 
    ${\bar f}^*=(f^*_i)_{i< k}$ a $k$-sequence of elements of $\omega\ho$, and
    $\n{\bar f}=(\n f_i)_{i< k}$ a $k$-sequence of $P_\beta$-names 
    of elements of $\myC$.\\
    We say ``${\bar f}^*$ is an honest 
    $(\alpha,\gamma,\beta)$-approximation of $\n{\bar f}$ 
    witnessed by $(p^m)_{m\in\omega}$'' if
    $(p^m)_{m\in\omega}$ is 
    an honest $(\alpha,\gamma,\beta)$-sequence
    and $p^m\forc_{(\alpha,\beta)}
    \n f_i\restriction m=f^*_i\restriction m$
    for all $m\in\omega$ and $i<k$.
  \item   
    ``${\bar f}^*$ is an honest $(\alpha,\beta)$-approximation of
    $\n{\bar f}$ under $p$'' means that there is a 
    $\gamma$ and a $(p^m)_{m\in\omega}$ such that $p^0\leq p$ and
    ${\bar f}^*$ is an honest     
    $(\alpha,\gamma,\beta)$-approximation of $\n{\bar f}$    
    witnessed by $(p^m)_{m\in\omega}$.
  \end{itemize}
\end{Def}

\begin{Lem}\label{lem:triv1} Let $\alpha\leq\zeta\leq\beta\leq\epsilon$. $P_\alpha$ forces:
  \begin{enumerate}
  \item\label{item:restr}
    If $(p^m)_{m\in\omega}$ is an honest $(\alpha,\gamma,\beta)$-sequence,
    then $(p^m\restriction\zeta)_{m\in\omega}$ 
    is an honest \linebreak[4]$(\alpha,\min(\zeta,\gamma),\zeta)$\nobreakdash-sequence.
  \item\label{item:thereexistsinterpr} Assume that 
      $p$               is an element of $P_\beta/G_\alpha$,
      $k$               a natural number,
      $(\n f_i)_{i< k}$ a $k$-sequence of $P_\beta$-names for
                        elements of $\myC$, and
      $D$               a dense subset of $P_\beta/G_\alpha$.
    Then there are $p'\leq p$ in $D$ 
    and $(f^*_i)_{i<k}$
    such that $(f^*_i)_{i<k}$ is an honest
    $(\alpha,\beta)$-approximation of $(\n f_i)_{i<k}$ under $p'$.
  \end{enumerate}
\end{Lem}

\begin{proof}
  We just show (\ref{item:thereexistsinterpr}). Work in $V_\alpha$.

  Let $\alpha\leq\gamma < \beta$ be minimal such that $p\restriction
  \gamma\not\forc_{(\alpha,\gamma)} p(\gamma)\in \an Q^\sigma_\gamma$.  If
  there is no such  $\gamma$, set $\gamma=\beta$ and $p_2=p$.
  Otherwise pick an $r\leq p\restriction \gamma$ in 
  $P_\gamma/G_\alpha$ such that
  $r\forc_{(\alpha,\gamma)} p(\gamma)\in \an Q^{\lnot\sigma}_\gamma$, and
  set $p_2=p\wedge r$. 

  Pick $p'\leq p_2$ in $D$.

  Let $\bar f^*$ approximate $\n{\bar f}$ witnessed by
  $p'=q^0\geq q^1\geq \dots$ (in $P_\beta/G_\alpha$).
  According to Lemma \ref{lem:quotient}(\ref{item:triv3252})
  there is a $q^\omega\in P_\gamma/G_\alpha$ such that
  $q^\omega\leq p'\restriction \gamma$ and
  $q^\omega\forc_{(\alpha,\gamma)} q^m\restriction\gamma\in G_\gamma$
  for all $m$.
  If $\gamma<\beta$, we can assume that 
  $q^\omega$ decides whether $ \{q^m(\gamma):\, m\in \omega\}$ is consistent.

  Set $r^m=q^\omega\wedge (q^m \restriction \beta\setminus \gamma)$, cf. \ref{lem:quotient}(\ref{item:modify1}).

  Assume $\gamma<\beta$ and $q^\omega$ forces consistency, i.e. $q^\omega\forc_{(\alpha,\gamma)}
  s\leq r^m(\gamma)$ for all $m$. Then $q^\omega$ forces that there is 
  an inconsistent sequence $s=s^0\geq s^1\geq\dots$ (since $s\in Q^{\lnot \sigma}_\gamma$).
  Modify $r^m$ such that $r^m\restriction \gamma=q^\omega\forc r^m(\gamma)=s^m$.
\end{proof}

\begin{ILem}\label{inductionlem}
  Assume that $q\in P_\alpha$ and that the following are in $N$:\\
  \myeinr $\alpha\leq\beta\leq \epsilon$, the $P_\alpha$-names $\an p$, $\an k$,
          $\an{\bar f}^*=(\an f_i^*)_{i\in \an k}$ and the $P_\beta$-name
          $\n{\bar f}=(\n f_i)_{i\in \an k}$ for elements of $\myC$.\\
  Assume that $q$ forces
  \begin{itemize}
  \item $\an{\bar f}^*$ is an honest
    $(\alpha,\beta)$-approximation of $\n{\bar f}$ under $\an p$ (in particular
    $\an p\in P_\beta/G_\alpha)$,
  \item $G_\alpha$ is $N$-generic and $\eta$ covers $N[G_\alpha]$.
  \end{itemize}
  Then there is a $q^+\in P_\beta$ such that $q^+\restriction \alpha=q$
  and $q^+$ forces
  \begin{itemize}
  \item $\an p\in G_\beta$,
  \item $G_\beta$ is $N$-generic and $\eta$ covers $N[G_\beta]$,
  \item $\an f^*_i\mR_j \eta$ implies $\n f_i \mR_{j} \eta$
        for all $i\in k$, $j\in\omega$.
  \end{itemize}
\end{ILem}

We prove the lemma by induction on~$\beta$. For $\alpha=\beta$
there is nothing to do. We split the proof into two cases:
$\beta$ successor and $\beta$ limit.

\begin{proof}[Proof for the case $\beta=\zeta+1$ successor]
  Let $\an p^m$ be $P_\alpha$-names for witnesses of the approximation.

  First assume that $q\in G_\zeta$ (i.e. $q\in G_\zeta\cap P_\alpha=G_\alpha$)
  and work in $V_\zeta$. 
  Set $p^{-1}=1_{P_\beta}$.  Let $-1\leq m*\leq
  \omega$ be
  the supremum of $\{m:\, \an p^m\restriction \zeta\in  G_\zeta\}$.

  Case 1: $m*=\omega$. In this case set ${\bar f}^{**}\DEFEQ \an{\bar f}^{*}$
        and $r\DEFEQ p^{0}(\zeta)\in \an Q_\zeta$. Note that 
        $\an p^m(\zeta)\forc_{\an Q_\zeta} \n f_i\restriction m = 
        f_i^{**}\restriction m$, i.e. ${\bar f}^{**}$ is an interpretation
        of $\bar{\n f}$ (with respect to $\an Q_\zeta$)
        under $r=\an p^{0}(\zeta)$.

  Case 2: $m*<\omega$. Find a $\an Q_\zeta$-interpretation ${\bar f}^{**}$
        of $\n {\bar f}$ under $r=\an p^{m*}(\zeta)\in\an Q_\zeta$ (use
        the fact the $\an Q_\zeta$ is preserving). Note that 
        $f_i^{**}\restriction m*=f_i^{*}\restriction m*$.

  Now fix (in $V$) $P_\zeta$-names $\an {\bar f}^{**}$ and $\an r$
  for this ${\bar f}^{**}$ and $r$ (we do not care how these names behave
  if $q\notin G_\alpha$). Then we get
  \[
    q\forc_\alpha \an p^m\restriction\zeta \forc_{(\alpha,\zeta)}
    \an f^{**}_i\restriction m=\an f^*_i\restriction m\text{ for all }i<\an k.
  \]
  So by fact \ref{lem:triv1}.(\ref{item:restr}), $q$ forces that
  $\an{\bar f}^*$ is an honest
  $(\alpha,\zeta)$-approximation of $\an{\bar f}^{**}$
  under $\an p\restriction\zeta$.

  By the induction hypothesis there is an $N$-generic $q^+\in P_\zeta$ which forces
  that $\an p^0\restriction\zeta\in G_\zeta$, $\eta$
  covers $N[G_\zeta]$ and of course that $\an Q_\zeta$ is proper and preserving.
  Assume $q^+\in G_\zeta$ and work in  $V_\zeta$. 
  Since $\an Q_\zeta$ is preserving and 
  ${\bar f}^{**}$ is an approximation of $\n{\bar f}$ under $r$,
  there is an $N[G_\zeta]$-generic 
  $q'\leq r$ which forces that $\eta$ covers $N[G_\zeta][G(\zeta)]$.
  Let (in $V$) $\an q'$ be a name for this $\an q$, and 
  set $q^{++}\DEFEQ q^+\wedge \an q'$. This $q^{++}$ is as required.
  (To see that $q^{++}\forc \an p \in G_\beta$, note that 
  $q^+\forc (\an p\restriction \zeta\in G_\zeta \,\&\, \an q'\leq \an p(\zeta))$.)
\end{proof}

\begin{proof}[Proof for the case $\beta$ limit]

  Choose a cofinal, increasing sequence $(\alpha_n)_{n\in\omega}$ 
  in $\beta\cap N$ such that $\alpha=\alpha_0$.

  Let $(D_n)_{n\in\omega}$ enumerate a basis of the open dense subsets
  of $P_\beta$ that are in $N$, and $(\n g_n)_{n\in\omega}$ all
  $P_\beta$-names in $N$ for elements of $\myC$. We may assume that
  $D_0=P_\beta$,
  $D_{n+1}\subseteq D_n$ and that every $p\in D_{n+1}$ decides
  $\n g_m\restriction n$ for $0\leq m\leq n$ as well as 
  $\an k$ and
  $\n f_i\restriction n$ for $0\leq i\leq \an k$.

  Let $\an\gamma_0$ and $(\an p^m_0)_{m\in\omega}$ be 
  $P_{\alpha_0}$-names for witnesses of the approximation
  in the assumption.
  Set $q_{-1}\DEFEQ q$, $\an k_0\DEFEQ \an k$
  and $\an {\bar f}^{*0}\DEFEQ \an {\bar f}^{*}$.

  Given $\an k_n$, we set $\n {\bar f}^{n}= (\n f^{n}_i)_{i<\an k_n} \DEFEQ(\n
  f_0,\dots,\n f_{\an k-1},\n g_0,\dots,\n g_{\an k_n-\an k})$.
  
  By induction on~$n\geq 1$ we can construct the following 
  $P_{\alpha_n}$-names in $N$:\\
  \begin{tabular}{ll}
    $(\an p^m_n)_{m\in\omega}$ & a sequence of conditions
                                 in $P_\beta/G_{\alpha_n}$,\\
    $\an \gamma_n$             & an ordinal,\\
    $\an k_n$                  & a natural number $\geq \an k_{n-1}$,\\
    $\an {\bar f}^{*n}=(\an f^{*n}_i)_{i<\an k_n}$
                               & a $\an k_n$-sequence of functions
                                 from $\omega$ to $\omega$,\\
  \end{tabular}\\
  such that  (for $n\geq 1$)
  $P_{\alpha_n}$ forces  that $\an p^0_{n-1}\restriction
  \alpha_n\in G_{\alpha_n}$ implies\footnote{or: 
  $\forc_{\alpha_{n-1}} p^0_{n-1}\restriction \alpha_n\forc_{(\alpha_{n-1},\alpha_n)}$}
  \begin{itemize}
    \item $\bar{\an f}^{*n}$ is an honest $(\alpha_n,\an\gamma_n,\beta)$-approximation
          of $\bar{\n f}^n$ witnessed by $(\an p^m_n)_{m\in\omega}$,
    \item One of the following cases holds:
      \begin{itemize}
      \item[$\an A_n$] $\an\gamma_{n-1}<\alpha_{n}$. Then
        there is a maximal $m*\geq 0$ such that $\an
        p^{m*}_{n-1}\restriction\alpha_n$ is in $G_{\alpha_n}$.
        Then we set $\an k_{n}\DEFEQ n+\an k$ and choose $\an p^0_n\leq \an
        p^{m*}_{n-1}\leq p^0_{n-1}$, $p^0_n\in D_n$.  
      \item[$\an B_n$] $\an\gamma_{n-1}= \beta$.
        (In this case the rest of the iteration is $\sigma$-complete and
        all $\an p^m_{n-1}$ are identical.)
        Set $\an k_{n}\DEFEQ n+\an k$ and choose 
        $\an p^0_n\leq p^0_{n-1}$ in $D_n$.
      \item[$\an C_n$] $\alpha_n\leq \an\gamma_{n-1}< \beta$. (Then all
        $\an p^m_{n-1}\restriction {\alpha_n}$ are identical and
        therefore in $P_{\alpha_n}/G_{\alpha_n}$.)
        In this case we ``do nothing'', i.e. we set
        $p^m_n\DEFEQ p^m_{n-1}$, $\an k_{n}\DEFEQ \an k_{n-1}$ and
        $\an {\bar f}^{*n}\DEFEQ \an {\bar f}^{*n-1}$.
      \end{itemize}
  \end{itemize}
  All we need for this construction is 
  \ref{lem:triv1}(\ref{item:thereexistsinterpr}).
  Note that in all three cases $\an p^0_n\leq \an p^0_{n-1}$;
  in case $\an A_n$ or $\an B_n$ $\an p^0_n\in D_n$ and therefore
  $\an p^0_n\forc_{(\alpha_n,\beta)} \n f^n_i\restriction n=\an 
  f^{\ast n}\restriction n$ for $i<n$. In case $\an B_n$, $\an \gamma_n$
  is again $\beta$, in case $\an C_n$, $\an \gamma_n=\an\gamma_{n-1}$.
  In all three cases, $f^{\ast n}$ is an honest 
  $(\alpha_n,\gamma_n,\alpha_{n+1})$-approximation witnessed by
  $(\an p^m_n\restriction\alpha_{n+1})_{m\in\omega}$.

  To see this, we just have to show that 
  $\an p^m_{n}\restriction \alpha_{n+1}
  \forc_{(\alpha_{n},\alpha_{n+1})}
  \an{f}_i^{*n+1}\restriction m=\an{f}_i^{*n}\restriction m$.
  Assume $G_{\alpha_{n+1}}$ contains 
  $p^m_{n}\restriction \alpha_{n+1}$. 
  Then in $V_{\alpha+2}$, 
  case $A_{n+1}$, $B_{n+1}$ or $C_{n+1}$ holds. In each case 
  we can 
  extend $G_{\alpha_{n+1}}$ to a $P_\beta$-generic filter
  $G_\beta$ containing $\an p^m_{n+1}$. Then (by case distinction)
  $G_\beta$ contains $\an p^m_{n}$ as well, i.e. 
  $\an{f}_i^{*n}\restriction m=\n f_i\restriction m=
  \an{f}_i^{*n+1}\restriction m$.
  
  Next we construct (by induction on~$n\geq 0$)
  $q_n\in P_{\alpha_{n+1}}$ such that $q_{n}\restriction \alpha_{n}=q_{n-1}$ 
  and $q_n$ forces:
  \begin{itemize}
  \item $G_{\alpha_{n+1}}$ is $N$-generic and $\eta$ covers
     $N[G_{\alpha_{n+1}}]$,
  \item 
    $\an p^0_{n}\restriction \alpha_{n+1}\in G_{\alpha_{n+1}}$,
  \item 
    $\an f^{*n}_i\mR_{j}\eta$ implies $\an f^{*n+1}_i \mR_{j} \eta$ for
    $i\in \an k_{n}$, $j\in\omega$,
  \item $(\an{f}_i^{*n+1})_{i<\an k_{n+1}}$ approximates
    $(\an{f}_i^{*n+2})_{i<\an k_{n+1}}$
    witnessed by $(\an p^m_{n+1}\restriction \alpha_{n+2})_{m\in\omega}$.
  \end{itemize}

  We can do this simply by applying the induction lemma iteratively:
  Given $q_{n-1}$,
  we choose $q_n$ using  \ref{inductionlem} as induction hypothesis,
  setting $\alpha\DEFEQ \alpha_n$, $\beta\DEFEQ \alpha_{n+1}$,
  $q\DEFEQ q_{n-1}$, $q^+\DEFEQ q_n$,
  $\an p\DEFEQ \an p^0_n$,
  $\an k\DEFEQ \an k_{n}$,
  $\an{\bar f}^*\DEFEQ \an{\bar f}^{*n}$,
  $\bar {\n f}\DEFEQ \an{\bar f}^{*n+1}$. 
  
  Now $q_\beta\DEFEQ \bigcup q_{\alpha_n}$ is as required:
  Assume $G_\beta$ is a $P_\beta$-generic filter over $V$
  containing $q_\beta$. We write $p^m_n$ for $\an p^m_n[G_\beta]=
  \an p^m_n[G_{\alpha_n}]$ etc.
  \begin{itemize}
  \item $p^0_n\in G_\beta$ for all $n$:

     $q_m\forc \an p^0_{m-1}\restriction \alpha_m\in G_{\alpha_m}$ for 
     all $m$.  Therefore $p^0_m\leq p^0_{m-1}$ for all $m$. 
     So for $m>n$, $q_m\forc \an p^0_{n}\restriction \alpha_m\in G_{\alpha_m}$.
     Therefore $p^0_n\restriction \alpha_m\in G_{\alpha_m}$ for all
     $m$, i.e. $p^0_n\in G_\beta$.

  \item $\gamma_n=\gamma_{n-1}$ unless $\gamma_{n-1}<\alpha_n$ 
    (i.e. case $A_n$ holds).
  \item $\bigcup_{n\in\omega} {k_n}=\omega$, and infinitely often 
     case $A_n$ or case $B_n$ holds:

     If $\gamma_m=\beta$ for some $m$, then case $B_n$ holds 
     (and $k_n=n$) for all $n>m$. Whenever $\alpha_{m+1}\leq \gamma_m<\beta$
     (i.e. case $C_{m+1}$ holds), then for some $n>m$
     (the smallest $n$  such that $\alpha_n>\gamma_m$)
     case $A_n$ holds and therefore $k_n=n$.

  \item $G_\beta$ is $N$-generic.

     Let $D\in N$ be dense. Then  $D\supseteq D_m\in N$, and
     for some $n\geq m$, case $A_n$ or case $B_n$ holds.
     Therefore $p^0_n\in N\cap D_n\cap G_\beta$, and $D_n\subseteq D_m$.

  \item We set $f^\infty_i\DEFEQ \an f^l_i[G_\beta]$ for some
     $l$ sufficiently large (i.e. $l$ such that $k_l>i$).

     So  $(f^\infty_0,f^\infty_1,\dots)=(f_0,\dots,f_{k-1},g_0,g_1,\dots)$.

  \item If $k_n>i$ and $l>n$, then $f^{*n}_i \mR_j \eta$ implies  $f^{*l}_i \mR_j \eta$.

  \item If $k_n>i$, then $f^{*n}\mR_j\eta$ implies $f^\infty_i\mR_j\eta$.

    Recall that $\{f:\, f \mR_j\eta\}$ is closed. For every $m$ there
    there is an $l>m$ such that case $A_l$ or $B_l$ holds, i.e. 
    $f^{*l}_i\restriction l=f^\infty_i\restriction l$, and by the last
    item $f^{*l}_i\mR_j\eta$.

  \item $\eta$ covers $N[G_\beta]$.

     Let $g\in N[G_\beta]\cap \myC$. Then for some $i$,
     $g=f^\infty_i$. Pick an $n$ such that $k_n>i$.  Since $\eta$ covers
     $N[G_{\alpha_n}]$ and $\an f^{*n}_i\in N[G_{\alpha_n}]$, $\an f^{*n}_i\mR_j
     \eta$ for some $j\in\omega$.
     
  \end{itemize}
  This ends the proof of the limit case.
\end{proof}

Note that the iteration lemma applied to the case $\alpha=0$ does not
immediately give the  preservation theorem \ref{thm:martin}, 
since we only get preservation for
honest approximations. This turns out to be no problem, however. 
Let us recall the structure of the proof:
  
Assume that $(P^0_\alpha,\an Q^0_\alpha)_{\alpha\in\epsilon}$
is a proper countable support iteration such that
$P^0_\alpha$ forces that $Q^0_\alpha$ is densely preserving
for all $\alpha$.
\begin{itemize}
\item Define $P^0_\alpha$-names $\an Q^1_\alpha$ 
  so that $P^0_\alpha$ forces that $Q^1_\alpha$ is a dense subforcing
  of $Q^0_\alpha$ and preserving (we can do that by the definition
  of densely preserving).
\item Define $P^0_\alpha$-names $Q^2_\alpha$ 
  so that $P^0_\alpha$ forces that $Q^2_\alpha$ is a dense subforcing
  of $Q^1_\alpha$ and decisive (we can do that by Fact \ref{fact:densesubpreservingagain}).
  $Q^2_\alpha$ is still preserving by Fact \ref{fact:decisiveisdense}. 
\item Let $(P_\alpha,\an Q_\alpha)$ be the countable support
  iteration as in Fact \ref{denseiswurscht}, obtained from $Q^2_\alpha$. 
  In particular $P_\alpha$ forces that $Q_\alpha$ is
  decisive and preserving (so we can apply the induction
  lemma), and 
  $P_\alpha$ can be densely embedded into $P^0_\alpha$
  for all $\alpha\leq \epsilon$.
\item Set 
  $P'\DEFEQ\{1_{P_\epsilon}\}\cup\{p\in P_\epsilon:\, 
  (\exists \gamma\leq \epsilon)\, 
  (\gamma=\epsilon\,\vee\, 
  p\restriction\gamma\forc_\gamma p(\gamma)\in Q^{\lnot\sigma}) 
  \,\&\,$\\
  \phantom{x}\hfill{$(\forall \alpha<\gamma) 
  p\restriction\alpha\forc_\alpha p(\alpha)\in Q^\sigma
  \}$.}\\
  $P'$ is a dense subforcing of $P_\epsilon$ and therefore of $P^0_\epsilon$.
  We assign to every $p\in P'\setminus\{1_{P'}\}$ the (unique) corresponding
  $\gamma(p)$. If $q\leq p$, then $\gamma(q)=\gamma(p)$.
\item 
  We claim that $P'$ is preserving (this finishes the proof
  of the iteration theorem). 
  Assume that (in $P'$)
  ${\bar f}^*$ interprets $\n{\bar f}$ witnessed by $(p^m)_{m\in\omega}$.
  We have to show that there is an honest witness
  $(p_1^m)_{m\in\omega}$ such that $p_1^0\leq p^0$.
  \begin{itemize}
  \item
    If all $p^m$ are $1_P$, then $\n{\bar f}$ is the standard name for
    ${\bar f}^*$ and there is nothing to do. So
    let $m*$ be the smallest $m$ such that $p^{m*}\neq 1_P$. Set
    $\gamma=\gamma(p^{m*})$.
  \item
    There is a $p^\omega$ in $P_\gamma$ such that
    $p^\omega\leq p^m\restriction \gamma$ for all $m$.
    Set $p^m_1\DEFEQ p^\omega \wedge p^m$. (So if 
    $\gamma=\epsilon$, then $p^m_1=p^\omega$ for all $m$.)
  \item
    If $\gamma<\epsilon$, we can assume that  $p^\omega$ decides
    whether the set $\{\an p^m(\gamma):\, m\in\omega\}$ is consistent. 
    If it decides positively, then we redefine
    $\an p_1^m(\gamma)$ to be any inconsistent sequence 
    in $\an Q_\gamma$ stronger than all $\an p^m(\gamma)$.
  \item 
    The resulting sequence $(p_1^m)_{m\in\omega}$
    witnesses that $\bar f^*$ is an honest approximation
    of $\n {\bar f}$.
  \end{itemize}
\end{itemize}

\bibliographystyle{plain}
\bibliography{tp}

\end{document}